\documentclass[11pt]{amsart}
\usepackage{amssymb,amsmath} 
\usepackage[all]{xy}
\title[Inverting weak dihomotopy equivalence]{Inverting weak dihomotopy equivalence using homotopy continuous flow}
\author[P. Gaucher]{Philippe Gaucher}
\address{Preuves Programmes et Syst{\`e}mes\\ Universit{\'e} Paris 7--Denis Diderot\\
Case 7014\\2 Place Jussieu\\ 75251 PARIS Cedex 05\\ France}
\email{gaucher@pps.jussieu.fr}
\urladdr{http://www.pps.jussieu.fr/{\~{}}gaucher/}

\newcommand{\C}{\mathcal{C}}
\newcommand{\p}\times
\renewcommand{\P}{\mathbb{P}}
\newcommand{\beas}{\begin{eqnarray*}}
\newcommand{\eeas}{\end{eqnarray*}}
\newtheorem{thm}{Theorem}
\newtheorem{prop}{Proposition}
\newtheorem{lem}{Lemma}
\newtheorem{cor}{Corollary}
\newtheorem{defn}{Definition}
\newtheorem{nota}{Notation}
\newcommand{\bd}{\begin{defn}}
\newcommand{\ed}{\end{defn}}
\newcommand{\bp}{\begin{prop}}
\newcommand{\ep}{\end{prop}}
\newcommand{\bth}{\begin{thm}}
\renewcommand{\eth}{\end{thm}}
\newcommand{\bpf}{\begin{proof}}
\newcommand{\epf}{\end{proof}}
\newcommand{\fr}[1]{\ar@{->}[r]^-{#1}}
\newcommand{\fd}[1]{\ar@{->}[d]_{#1}}
\renewcommand{\top}{{\mathbf{Top}}}
\newcommand{\ho}{{\mathbf{Ho}}}
\newcommand{\iso}{\cong}
\newcommand{\T}{{\mathcal{T}}}
\newcommand{\brm}[1]{\rm{\mathbf{#1}}}
\newcommand{\dtop}{{\brm{Flow}}}
\newcommand{\glob}{{\rm{Glob}}}
\newcommand{\liminj}{\varinjlim}

\def\cocartesien{%
  \ar@{-}[]+L+<-6pt,+2pt>;[]+LU+<-6pt,+6pt>%
  \ar@{-}[]+U+<-2pt,+6pt>;[]+LU+<-6pt,+6pt>%
}

\DeclareMathOperator{\dom}{dom}
\DeclareMathOperator{\codom}{codom}
\DeclareMathOperator{\id}{Id}
\DeclareMathOperator{\cell}{{\brm{cell}}}
\DeclareMathOperator{\cof}{{\brm{cof}}}
\DeclareMathOperator{\inj}{{\brm{inj}}}
\DeclareMathOperator{\Path}{{\rm{Path}}}
\DeclareMathOperator{\Iso}{{\underline{Iso}}}
\DeclareMathOperator{\All}{{\underline{All}}}

\begin{document}

\begin{abstract} 
  A flow is homotopy continuous if it is indefinitely divisible up to
  S-homo\-to\-py. The full subcategory of cofibrant homotopy continuous
  flows has nice features. Not only it is big enough to contain all
  dihomotopy types, but also a morphism between them is a weak
  dihomotopy equivalence if and only if it is invertible up to
  dihomotopy. Thus, the category of cofibrant homotopy continuous
  flows provides an implementation of Whitehead's theorem for the full
  dihomotopy relation, and not only for S-homotopy as in previous
  works of the author. This fact is not the consequence of the
  existence of a model structure on the category of flows because it
  is known that there does not exist any model structure on it whose
  weak equivalences are exactly the weak dihomotopy equivalences. This
  fact is an application of a general result for the localization of a
  model category with respect to a weak factorization system.
\end{abstract}

\maketitle
\tableofcontents

\section{Introduction}

There are numerous uses of the notion of ``category without
identities''. For recent papers, see for example \cite{MR1769341}
\cite{MR1928230} \cite{MR2153892}. An enriched version of this notion,
in the sense of \cite{MR651714}, over the category of general
topological spaces can be found in \cite{qalg9608025}. By considering
``small categories without identities'' enriched over the category of
compactly generated topological spaces, that is weak Hausdorff
k-spaces in the sense of \cite{MR99h:55031}, one obtains an object
called a \textit{flow} which allows a model categorical treatment of
\textit{dihomotopy} (directed homotopy). Indeed, a flow $X$ can model
(the time flow of) a higher dimensional automaton \cite{Pratt}
\cite{rvg} \cite{exHDA} as follows. A flow $X$ consists of
\begin{enumerate} 
\item a set of \textit{states} $X^0$; 
\item for each pair of states $(\alpha,\beta)\in X^0 \p X^0$, there is
  a compactly generated topological space $\P_{\alpha,\beta}X$ called
  the \textit{path space between $\alpha$ and $\beta$} representing
  the concurrency between $\alpha$ and $\beta$; each element of
  $\P_{\alpha,\beta}X$ corresponds to a \textit{non-constant execution
    path} from $\alpha$ to $\beta$; the
  emptiness of the space $\P_{\alpha,\alpha}X$ for some state $\alpha$
  means that there are no loops from $\alpha$ to itself; let
\[\P X=  \bigsqcup_{(\alpha,\beta)\in X^0\p X^0}\P_{\alpha,\beta}X.\] 
\item for each triple of states $(\alpha,\beta,\gamma)\in X^0 \p X^0
  \p X^0$, there is a strictly associative composition law
  $\P_{\alpha,\beta}X \p \P_{\beta,\gamma}X \longrightarrow
  \P_{\alpha,\gamma}X$ corresponding to the concatenation of
  non-constant execution paths.
\end{enumerate}

\begin{figure}
\[
\xymatrix{ \widehat{0} \ar@/^20pt/[r] \ar@/_20pt/[r] &  \widehat{1}\\
\widehat{0} \ar@(ur,dr) &
}
\]
\caption{Bad identification of a $1$-dimensional empty globe and a loop after 
a contraction in the direction of time}
\label{ex2}
\end{figure}
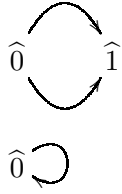

The main problem to model dihomotopy is that contractions in the
direction of time are forbidden. Otherwise in the categorical
localization of flows with respect to the dihomotopy equivalences, the
relevant geometric information is lost \cite{1eme} \cite{diCW}.  Here
is a very simple example. Take two non-constant execution paths going
from one initial state to one final state.  If contractions in the
direction of time were allowed, then one would find in the same
equivalence class a loop (cf. Figure~\ref{ex2}): this is not
acceptable.

Two kinds of deformations are of interest in the framework of flows.
The first one is called \textit{weak S-homotopy equivalence}: it is a
morphism of flows $f:X\longrightarrow Y$ such that the set map
$f^0:X^0\longrightarrow Y^0$ is a bijection and such that the
continuous map $\P f:\P X \longrightarrow \P Y$ is a weak homotopy
equivalence. It turns out that there exists a model structure on the
category of flows whose weak equivalences are exactly the weak
S-homotopy equivalences (\cite{model3} and Section~\ref{hc} of this
paper). However, the identifications allowed by the weak S-homotopy
equivalences are too rigid. So another kind of weak equivalence is
required. The \textit{T-homotopy equivalences} are generated by a set
$\mathcal{T}$ of cofibrations obtained by taking the cofibrant
replacement of the inclusions of posets~\footnote{Any poset $P$ can be
  viewed as a flow in an obvious way: the set of states is the
  underlying set of $P$ and there is a non-constant execution path
  from $\alpha$ to $\beta$ if and only if $\alpha<\beta$. Note that
  the inequality is strict. Indeed, the ordering of $P$ represents the
  direction of time and the flow associated with a poset must be
  loopless.} of Definition~\ref{genT}.  This approach of T-homotopy is
presented for the first time in \cite{1eme}. The latter models
``refinement of observation''. For instance, the inclusion of posets
$\{\widehat{0} < \widehat{1}\} \subset \{\widehat{0} < A <
\widehat{1}\}$ corresponds to the identification of a directed segment
$U$ going from the initial state $\widehat{0}$ to the final state
$\widehat{1}$ with the composite $U'*U''$ of two directed segments
(cf.  Figure~\ref{ex1}).

\begin{figure}
\[
\xymatrix{\widehat{0} \ar@{->}[rrrr]^-{U} &&&& \widehat{1} \\
\widehat{0} \ar@{->}[rr]^-{U'} && A \ar@{->}[rr]^-{U''} && \widehat{1}}
\]
\caption{The simplest example of refinement of observation}
\label{ex1}
\end{figure}

The problem we face can then be presented as follows. We have:
\begin{enumerate}
\item A model structure on the category of flows $\dtop$, called the
  weak S-homotopy model structure, such that the class of weak
  equivalences is exactly the class $\mathcal{S}$ of weak S-homotopy
  equivalences.  One wants to invert the weak S-homotopy equivalences
  because two weakly S-homotopy equivalent flows are equivalent from
  an observational viewpoint. This model structure provides an
  implementation of Whitehead's theorem for S-homotopy only.
\item A set of cofibrations $\T$ of \textit{generating T-homotopy
    equivalences} one would like to invert because these maps model
  refinement of observation.
\item Three known invariants with respect to weak S-homotopy and
  T-homotopy: the \textit{underlying homotopy type functor}
  \cite{model2}, the \textit{branching homology} and the
  \textit{merging homology} \cite{exbranch}.
\item \label{4} Every model structure on $\dtop$ which contains as
  weak equivalences the class of morphisms $\mathcal{S}\cup \T$,
  contains weak equivalences which do not preserve the three known
  invariants \cite{nonexistence}. In particular, the category
  $\dtop[\mathcal{S}_\T^{-1}]$ below is not the Quillen homotopy
  category of a model structure of $\dtop$. The left Bousfield
  localization of the weak S-homotopy model structure with respect to
  the set of cofibrations $\T$ is therefore not relevant here.
\end{enumerate}
The negative result (\ref{4}) prevents us from using the machinery of
model category on the category $\dtop$ for understanding the full
dihomotopy equivalence relation. There are then several possibilities:
reconstructing some pieces of homotopy theory in the framework of
flows, finding new categories for studying S-homotopy and T-homotopy,
or also relating dihomotopy on $\dtop$ to other axiomatic
presentations of homotopy theory. The possibility which is explored in
this paper is the first one.

Indeed, the goal of this work is to prove that it is possible to find
a full subcategory of the category of flows which is big enough to
contain all dihomotopy types and in which the weak dihomotopy
equivalences are exactly the invertible morphisms up to dihomotopy.
The main theorem of the paper states as follows (cf. Section~\ref{hc}
for a reminder about flows):

\begin{thm} (Theorem~\ref{fin3} and Theorem~\ref{fin4}) Let $J^{gl}$
  be the set of generating trivial cofibrations of the weak S-homotopy
  model structure of $\dtop$. Let $\T$ be the set of generating
  T-homotopy equivalences. Let $\dtop_{cof}$ be the full subcategory
  of cofibrant flows. There exists a full subcategory
  $\dtop_{cof}^{f,\T}$ of the category of cofibrant flows
  $\dtop_{cof}$, the one of homotopy continuous flows, a class of
  morphisms of flows $\mathcal{S}_\T$ and a congruence $\sim_\T$ on
  the morphisms of $\dtop$ such that the inclusion functors
  $\dtop_{cof}^{f,\T} \subset \dtop_{cof} \subset \dtop$ induce the
  equivalences of categories
\[\dtop_{cof}^{f,\T}/\!\sim_\T \simeq \dtop_{cof}[(\mathcal{S}\cup
\cof(J^{gl}\cup T))^{-1}] \simeq \dtop[\mathcal{S}_\T^{-1}].\]
Moreover, one has: 
\begin{enumerate}
\item The class of morphisms $\mathcal{S}_\T$ contains the weak
  S-homotopy equivalences and the morphisms of $\cof(J^{gl}\cup T)$
  with cofibrant domains.
\item Every morphism of $\mathcal{S}_\T$ preserves the underlying
  homotopy type, the branching homology and the merging homology.
\end{enumerate}
\end{thm}

We now outline the contents of the paper.  The purpose of
Section~\ref{ab} is to give the proof of the theorem above in a more
abstract setting. The starting point is a model category $\mathcal{M}$
together with a weak factorization system $(\mathcal{L}, \mathcal{R})$
satisfying some technical conditions which are fulfilled by the weak
S-homotopy model structure of $\dtop$ and by the set of generating
T-homotopy equivalences.  Several proofs of Section~\ref{ab} are
adaptations of standard proofs \cite{MR36:6480} \cite{MR99h:55031}.
But since the existence of a convenient model structure for
$(\mathcal{L},\mathcal{R})$ is not supposed~\footnote{that is: a model
  structure such that $\mathcal{L}$ is the class of trivial
  cofibrations.}, there are some subtle differences and also new
phenomena. The idea of considering the path object construction comes
from the reading of Kurz and Rosick\'y's paper \cite{ideeloc}. In this
paper, Kurz and Rosick\'y have the idea of considering a cylinder
object construction with any weak factorization system
$(\mathcal{L},\mathcal{R})$. This allows them to investigate the
categorical localization of the underlying category with respect to
the class of morphisms $\mathcal{R}$ viewed, morally speaking, as a
class of trivial fibrations. The dual situation is explored in this
section, with an underlying category which is not only a category but
also a model category. The situation described in Section~\ref{ab}
makes one think of the notion of fibration category in the sense of
Baues \cite{homotopieabstraite0}. However, we do not know how to
construct a fibration category from the results of Section~\ref{ab}.
The path object functor constructed in Section~\ref{ab} cannot satisfy
the whole set of axioms of a P-category in the sense of Baues
\cite{homotopieabstraite0} since the associated homotopy relation is
not transitive. In particular, it does not even seem to satisfy the
pullback axiom.  Next, Section~\ref{hc} proves the theorem above as an
application of Section~\ref{ab}.

\subsubsection*{Link with the series of papers ``T-homotopy and refinement
of observation''}

This paper is independent from the series of papers ``T-homotopy and
refinement of observation'' except for the proof of Theorem~\ref{fin4}
at the very end of this work in which \cite{1eme} Theorem~5.2 is used.
This paper was written while the author was trying to understand
whether the (categorical) localization $\dtop[\cof(\T)^{-1}]$ of the
category of flows with respect to the T-homotopy equivalences
introduced in \cite{1eme} is locally small. Indeed, the local
smallness is not established in the series of papers ``T-homotopy and
refinement of observation''. The result we obtain is more subtle. In
the ``correct'' localization, all morphisms of $\cof(J^{gl}\cup \T)$
with cofibrant domains are inverted. This is enough for future
application in computer science since the real concrete examples are
all of them modelled by cofibrant flows. But it is not known whether
the other morphisms of $\cof(J^{gl}\cup \T)$ are inverted. If this
fact should be true, then it would probably be a consequence of the
left properness of the weak S-homotopy model structure of $\dtop$
(which is proved in \cite{2eme} Theorem~6.4).

\section{Prerequisites and notations}

The initial object (resp. the terminal object) of a category $\C$, if
it exists, is denoted by $\varnothing$ (resp. $\mathbf{1}$).

Let $i:A\longrightarrow B$ and $p:X\longrightarrow Y$ be maps in a
category $\C$. Then $i$ has the {\rm left lifting property} (LLP) with
respect to $p$ (or $p$ has the {\rm right lifting property} (RLP) with
respect to $i$) if for every commutative square
\[
\xymatrix{
A\ar@{->}[dd]_{i} \ar@{->}[rr]^{\alpha} && X \ar@{->}[dd]^{p} \\
&&\\
B \ar@{-->}[rruu]^{g}\ar@{->}[rr]_{\beta} && Y,}
\]
there exists a morphism $g$ called a \textit{lift} making both
triangles commutative.

Let $\C$ be a cocomplete category.  If $K$ is a set of morphisms of
$\C$, then the class of morphisms of $\C$ that satisfy the RLP
(\textit{right lifting property}) with respect to every morphism of
$K$ is denoted by $\inj(K)$ and the class of morphisms of $\C$ that
are transfinite compositions of pushouts of elements of $K$ is denoted
by $\cell(K)$. Denote by $\cof(K)$ the class of morphisms of $\C$ that
satisfy the LLP (\textit{left lifting property}) with respect to every
morphism of $\inj(K)$.  The cocompleteness of $\C$ implies
$\cell(K)\subset \cof(K)$. Moreover, every morphism of $\cof(K)$ is a
retract of a morphism of $\cell(K)$ as soon as the domains of $K$ are
small relative to $\cell(K)$ (\cite{MR99h:55031} Corollary~2.1.15). An
element of $\cell(K)$ is called a \textit{relative $K$-cell complex}.
If $X$ is an object of $\C$, and if the canonical morphism
$\varnothing\longrightarrow X$ is a relative $K$-cell complex, one
says that $X$ is a \textit{$K$-cell complex}.

A \textit{congruence} $\sim$ on a category $\C$ consists of an
equivalence relation on the set $\C(X,Y)$ of morphisms from $X$ to $Y$
for every object $X$ and $Y$ of $\C$ such that if $f,g\in\C(X,Y)$,
then $f\sim g$ implies $u\circ f\sim u\circ f$ and $f\circ v\sim
g\circ v$ for any morphism $u$ and $v$ as soon as $u\circ f$ and
$f\circ v$ exist.

Let $\C$ be a cocomplete category with a distinguished set of
morphisms $I$. Then let $\cell(\C,I)$ be the full subcategory of $\C$
consisting of the object $X$ of $\C$ such that the canonical morphism
$\varnothing\longrightarrow X$ is an object of $\cell(I)$. In other
terms, $\cell(\C,I)=(\varnothing\!\downarrow \! \C) \cap \cell(I)$.

It is obviously impossible to read this paper without some familiarity
with \textit{model categories}. Possible references for model
categories are \cite{MR99h:55031}, \cite{ref_model2} and
\cite{MR1361887}.  The original reference is \cite{MR36:6480} but
Quillen's axiomatization is not used in this paper. The Hovey's book
axiomatization is preferred.  If $\mathcal{M}$ is a
\textit{cofibrantly generated} model category with set of generating
cofibrations $I$, let $\cell(\mathcal{M}) := \cell(\mathcal{M},I)$.
Any cofibrantly generated model structure $\mathcal{M}$ comes with a
\textit{cofibrant replacement functor} $Q:\mathcal{M} \longrightarrow
\cell(\mathcal{M})$. For every morphism $f$ of $\mathcal{M}$, the
morphism $Q(f)$ is a cofibration, and even an inclusion of
subcomplexes. A set $K$ of morphisms of a model category
\textit{permits the small object argument} if the domains of the
morphisms of $K$ are small relative to $\cell(K)$. For such a set $K$,
one can use the small object argument. The small object argument is
recalled in the proof of Proposition~\ref{inclusionOK}.

In this paper, the notation $\simeq$ means \textit{weak equivalence}
or \textit{equivalence of categories}, and the notation $\iso$ means
\textit{isomorphism}.

A \textit{partially ordered set} $(P,\leq)$ (or \textit{poset}) is a
set equipped with a reflexive antisymmetric and transitive binary
relation $\leq$. A poset $(P,\leq)$ is \textit{bounded} if there exist
$\widehat{0}\in P$ and $\widehat{1}\in P$ such that $P\subset
[\widehat{0},\widehat{1}]$ and such that $\widehat{0} \neq
\widehat{1}$. Let $\widehat{0}=\min P$ (the bottom element) and
$\widehat{1}=\max P$ (the top element).

Every poset $P$, and in particular every ordinal, can be viewed as a
small category denoted in the same way: the objects are the elements
of $P$ and there exists a morphism from $x$ to $y$ if and only if
$x\leq y$. If $\lambda$ is an ordinal, a \textit{$\lambda$-sequence}
(or a \textit{transfinite sequence}) in a cocomplete category $\C$ is
a colimit-preserving functor $X$ from $\lambda$ to $\C$. We denote by
$X_\lambda$ the colimit $\liminj X$ and the morphism
$X_0\longrightarrow X_\lambda$ is called the \textit{transfinite
  composition} of the $X_\mu\longrightarrow X_{\mu+1}$.

If $\C$ is a locally small category, and if $\Sigma$ is a class of
morphisms of $\C$, then we denote by $\C[\Sigma^{-1}]$ the \textit{
  (categorical) localization} of $\C$ with respect to $\Sigma$
\cite{gz} \cite{MR1712872} \cite{MR96g:18001a}. The category
$\C[\Sigma^{-1}]$ is not necessarily locally small. If $\mathcal{M}$
is a model category with class of weak equivalences $\mathcal{W}$,
then the localization $\mathcal{M}[\mathcal{W}^{-1}]$ is locally small
and it is called the \textit{Quillen homotopy category} of
$\mathcal{M}$. It is denoted by $\ho(\mathcal{M})$.

\section{Localizing a model category w.r.t. a weak factorization system}
\label{ab}

\bd \cite{MR2003h:18001} Let $\C$ be a category. A {\rm weak
  factorization system} is a pair $(\mathcal{L},\mathcal{R})$ of
classes of morphisms of $\C$ such that the class $\mathcal{L}$ is the
class of morphisms having the LLP with respect to $\mathcal{R}$, such
that the class $\mathcal{R}$ is the class of morphisms having the RLP
with respect to $\mathcal{L}$ and such that every morphism of $\C$
factors as a composite $r\circ \ell$ with $\ell\in \mathcal{L}$ and
$r\in \mathcal{R}$. The weak factorization system is {\rm functorial}
if the factorization $r\circ \ell$ can be made functorial.  \ed

In a weak factorization system $(\mathcal{L},\mathcal{R})$, the class
$\mathcal{L}$ (resp.  $\mathcal{R}$) is completely determined by
$\mathcal{R}$ (resp.  $\mathcal{L}$).

\bd \label{cg} Let $\C$ be a cocomplete category. A weak factorization
system $(\mathcal{L},\mathcal{R})$ is {\rm cofibrantly generated} if
there exists a set $K$ of morphisms of $\C$ permitting the small
object argument such that $\mathcal{L}=\cof(K)$ and
$\mathcal{R}=\inj(K)$.  \ed

A cofibrantly generated weak factorization system is necessarily
functorial.  Definition~\ref{cg} appears in \cite{MR1780498} in the
context of locally presentable category as the notion of
\textit{small} weak factorization system.

\textit{\underline{The data for this section are}: 
\begin{enumerate}
\item a complete and cocomplete category $\mathcal{M}$ equipped with a
  model structure denoted by $({\rm Cof},{\rm Fib},\mathcal{W})$ for
  respectively the class of cofibrations, of fibrations and of weak
  equivalences such that the weak factorization system $({\rm Cof}\cap
  \mathcal{W},{\rm Fib})$ is cofibrantly generated : the set of
  generating trivial cofibrations is denoted by $J$.
\item a cofibrantly generated weak factorization system
  $(\mathcal{L},\mathcal{R})$ on $\mathcal{M}$ satisfying the
  following property: ${\rm Cof}\cap \mathcal{W}\subset \mathcal{L}
  \subset {\rm Cof}$. So there exists a set of morphisms $K$ such that
  $\mathcal{L} = \cof(J\cup K)$ and $\mathcal{R} = \inj(J\cup K)$ and
  such that $J\cup K$ permits the small object argument. Therefore
  every morphism $f$ factors as a composite $f = \beta(f) \circ
  \alpha(f)$ where $\alpha(f) \in \cell(J\cup K)\subset \mathcal{L}$
  and where $\beta(f) \in \mathcal{R}$. The functorial factorization
  is supposed to be obtained using the small object argument. It is
  fixed for the whole section.
\end{enumerate}}

\bd \label{pathdef} Let $X$ be an object of $\mathcal{M}$. The {\rm
  path object of $X$ with respect to $\mathcal{L}$} is the functorial
factorization
\[
\xymatrix{
X \ar@{->}[rr]^-{\alpha(\id_X,\id_X)}&& \Path_{\mathcal{L}}(X) \ar@{->}[rr]^-{\beta(\id_X,\id_X)} && X \p X}
\]
of the diagonal morphism $(\id_X,\id_X):X \longrightarrow X \p X$ by
the morphism $\alpha(\id_X,\id_X):X \longrightarrow
\Path_{\mathcal{L}}(X)$ of $\mathcal{L}$ composed with the morphism
$\beta(\id_X,\id_X):\Path_{\mathcal{L}}(X) \longrightarrow X \p X$ of
$\mathcal{R}$.  \ed

\begin{nota} 
  Let $\mathcal{M}_{cof}$ be the full subcategory of cofibrant objects
  of $\mathcal{M}$.
\end{nota}

The path object of $X$ with respect to $\mathcal{L}$ is cofibrant as
soon as $X$ is cofibrant since the morphism $\alpha(\id_X,\id_X): X
\longrightarrow \Path_{\mathcal{L}}(X)$ is a cofibration. So the path
object construction yields an endofunctor of $\mathcal{M}_{cof}$.

The reader must notice that we do not assume here that
$\alpha(\id_X,\id_X)$ is a weak equivalence of any kind, contrary to
the usual definition of a path object. As in \cite{ideeloc} for the
construction of the cylinder functor, we do use the functorial
factorization and we do suppose that $\alpha(\id_X,\id_X)$ belongs to
$\mathcal{L}$. A morphism of $\mathcal{L}$ being an isomorphism of
$\mathcal{M}_{cof}[(\mathcal{W} \cup \mathcal{L})^{-1}]$, our
condition is stronger than the usual one for the construction of a
path object in a model category.

\bd An object $X$ of $\mathcal{M}$ is {\rm fibrant with respect to
  $\mathcal{L}$} if the unique morphism $f_X:X \longrightarrow
\mathbf{1}$, where $\mathbf{1}$ is the terminal object of
$\mathcal{M}$, is an element of $\mathcal{R}$.  \ed

An object which is fibrant with respect to ${\rm Cof} \cap
\mathcal{W}$ is a fibrant object in the usual sense.

\begin{nota}
  Let $\mathcal{M}^{f,\mathcal{L}}$ be the full subcategory of
  $\mathcal{M}$ of fibrant objects with respect to $\mathcal{L}$.  Let
  $\mathcal{M}^{f,\mathcal{L}}_{cof}$ be the full subcategory of
  $\mathcal{M}_{cof}$ of fibrant objects with respect to
  $\mathcal{L}$.
\end{nota}

If $X$ is fibrant with respect to $\mathcal{L}$, the morphism $X\p X
\longrightarrow X \p \mathbf{1}\iso X$ belongs to $\mathcal{R}$.
Therefore the composite
\[\Path_\mathcal{L}(X)\longrightarrow X\p X \longrightarrow X \longrightarrow \mathbf{1}\] 
belongs to $\mathcal{R}$ as well. So the path object
$\Path_\mathcal{L}(X)$ is also fibrant with respect to
$\mathcal{L}$. Thus, the path object construction yields endofunctors
of $\mathcal{M}^{f,\mathcal{L}}$ and of
$\mathcal{M}^{f,\mathcal{L}}_{cof}$.

If $f:X \longrightarrow Y$ is a morphism of
$\mathcal{M}^{f,\mathcal{L}}_{cof}$, then the functorial factorization
$(\alpha,\beta)$ yields a composite a priori in $\mathcal{M}_{cof}$ 
(since $\mathcal{L} \subset {\rm Cof}$) 
\[
\xymatrix{
X\ar@{->}[rr]^-{\alpha(f)} && Z \ar@{->}[rr]^-{\beta(f)} && Y}
\]
equal to $f$. The unique morphism $Z \longrightarrow \mathbf{1}$ is
equal to the composite $Z\longrightarrow Y \longrightarrow \mathbf{1}$
of two morphisms of $\mathcal{R}$. Therefore $Z$ is fibrant with
respect to $\mathcal{L}$ and the functorial weak factorization system
$(\mathcal{L},\mathcal{R})$ restricts to a functorial weak
factorization system of $\mathcal{M}^{f,\mathcal{L}}_{cof}$ denoted in
the same way.

\bd \label{rightdef} Let $f,g:X \rightrightarrows Y$ be two morphisms
of $\mathcal{M}$. A {\rm right homotopy with respect to} $\mathcal{L}$
from $f$ to $g$ is a morphism $H: X \longrightarrow
\Path_{\mathcal{L}}(Y)$ such that
\[\beta(\id_Y,\id_Y)\circ H = (f,g).\] This situation is denoted by $f
\sim^r_{\mathcal{L}} g$.
\ed

Note the binary relation $\sim^r_{\mathcal{L}}$ does not depend on the
choice of the functorial factorization $(\alpha,\beta)$. Indeed, with
another functorial factorization $(\alpha',\beta')$, and the
corresponding path object functor $\Path'_\mathcal{L}$, one can
consider for every object $Y$ of $\mathcal{M}$ the commutative diagram
\[
\xymatrix{
Y \fd{} \ar@{=}[rr]{} && Y \fd{} \\
\Path_\mathcal{L}(Y) \fd{} \ar@{-->}[rr]^-{k}&&  \Path'_\mathcal{L}(Y) \fd{} \\
Y\p Y \ar@{->}[rr]^{\id_Y\p \id_Y} && Y\p Y.}
\] 
The lift $k$ exists since the arrow $Y \longrightarrow
\Path_\mathcal{L}(Y)$ is in $\mathcal{L}$ and since the arrow
$\Path'_\mathcal{L}(Y) \longrightarrow Y\p Y$ is in $\mathcal{R}$.

The morphism $\alpha(\id_Y,\id_Y)\circ f:X \longrightarrow
\Path_{\mathcal{L}} Y$ yields a right homotopy from $f$ to $f$ with
respect to $\mathcal{L}$. If $H:X \longrightarrow
\Path_{\mathcal{L}}(Y)$ is a right homotopy from $f$ to $g$ with
respect to $\mathcal{L}$, then the usual way for obtaining a right
homotopy from $g$ to $f$ with respect to $\mathcal{L}$ consists of
considering the commutative diagram:
\[
\xymatrix{
Y \fd{} \ar@{=}[rr]{} && Y \fd{} \\
\Path_\mathcal{L}(Y) \fd{} \ar@{-->}[rr]^-{k}&&  \Path_\mathcal{L}(Y) \fd{} \\
Y\p Y \ar@{->}[rr]^{\tau} && Y\p Y}
\] 
with $\tau(y,y')=(y',y)$. The existence of the lift $k$ comes from the
definition of the path object and of the fact that
$(\mathcal{L},\mathcal{R})$ is a weak factorization system. So the
binary relation $\sim^r_{\mathcal{L}}$ is reflexive and symmetric.

This relation is not transitive in general. The pair
$(\mathcal{R}^{op},\mathcal{L}^{op})$ is a weak factorization system
of the opposite category $\mathcal{M}^{op}$ (the model structure of
$\mathcal{M}$ is forgotten for this paragraph only). The path object
becomes a cylinder object and the binary relation $\sim^r_\mathcal{L}$
becomes the homotopy relation of \cite{ideeloc}. \cite{ideeloc}
Example~3.6 gives an example where the homotopy is not transitive.
Thus, the opposite category with the opposite weak factorization
system gives an example where $\sim^r_\mathcal{L}$ is not transitive.

\begin{nota} 
  Let us denote by $\sim_{\mathcal{L}}$ the transitive closure of the
  binary relation $\sim^r_{\mathcal{L}}$.
\end{nota}

\bp \label{extright} Let $X$ be an object of $\mathcal{M}_{cof}$. Let
$Y$ be an object of $\mathcal{M}$. Let $f,g:X\rightrightarrows Y$ be
two morphisms between them. Then $f\sim_{{\rm Cof}\cap \mathcal{W}} g$
if and only if $f$ and $g$ are right homotopic in the usual sense of
model categories.  \ep

Notice that it is crucial in the proof for $X$ to be cofibrant.

\bpf Indeed, two morphisms $f,g:X\rightrightarrows Y$ with $X$
cofibrant are right homotopic in the usual sense if the pair $(f,g)$
is in the transitive closure of the following situation denoted by
$f\sim^r g$ (cf. \cite{MR99h:55031} p7):
\begin{enumerate}
\item Decompose the diagonal morphism $(\id_Y,\id_Y):Y\longrightarrow
  Y\p Y$ into a weak equivalence $Y\longrightarrow PY$ of
  $\mathcal{M}$ followed by a fibration $(p_1,p_2):PY \longrightarrow
  Y\p Y$ of $\mathcal{M}$.
\item There exists $H:X \longrightarrow PY$ such that $(p_1, p_2)\circ
  H=(f,g)$.
\end{enumerate}
Let us factor the weak equivalence $Y \longrightarrow PY$ as a
composite $Y \longrightarrow P'Y \longrightarrow PY$ where $Y
\longrightarrow P'Y$ is a trivial cofibration and where $P'Y
\longrightarrow PY$ is a trivial fibration. Then one can lift the
right homotopy $H:X \longrightarrow PY$ to a morphism $\overline{H}:X
\longrightarrow P'Y$ since $X$ is cofibrant. But $\overline{H}$ is not
yet a right homotopy from $f$ to $g$ with respect to ${\rm Cof}\cap
\mathcal{W}$ since $P'Y$ is not necessarily the functorial path object
$\Path_{{\rm Cof}\cap \mathcal{W}}(Y)$ ! Let us consider the
commutative diagram
\[
\xymatrix{
Y \fd{} \ar@{=}[rr]{} && Y \fd{} \\
P'Y \fd{} \ar@{-->}[rr]^-{k}&&  \Path_{{\rm Cof}\cap \mathcal{W}}(Y) \fd{} \\
Y\p Y \ar@{->}[rr]^{\id_Y\p \id_Y} && Y\p Y}
\] 
Since the arrow $Y \longrightarrow P'Y$ is a trivial cofibration and
since the arrow $\Path_{{\rm Cof}\cap \mathcal{W}}(Y) \longrightarrow
Y\p Y$ is a fibration, there exists a lift $k$. Then $k\circ
\overline{H}$ is a right homotopy with respect to ${\rm Cof}\cap
\mathcal{W}$ from $f$ to $g$.

Conversely, the path object with respect to ${\rm Cof}\cap
\mathcal{W}$ is a path object in the above sense of model categories.
So a right homotopy from $f$ to $g$ with respect to ${\rm Cof}\cap
\mathcal{W}$ is a right homotopy in the usual sense of model
categories.  \epf

The following proposition gives a sufficient condition for the binary
relation $\sim^r_\mathcal{L}$ to be transitive.

\bp \label{traditionnel1} Let us suppose that there exists a model
structure $({\rm Cof}_\mathcal{L},{\rm
  Fib}_\mathcal{L},\mathcal{W}_\mathcal{L})$ on $\mathcal{M}$ such
that $\mathcal{L} = {\rm Cof}_\mathcal{L}\cap\mathcal{W}_\mathcal{L}$
and such that every cofibrant object of $\mathcal{M}$ is a cofibrant
object of $({\rm Cof}_\mathcal{L},{\rm
  Fib}_\mathcal{L},\mathcal{W}_\mathcal{L})$. Let $X$ and $Y$ be two
objects of $\mathcal{M}^{f,\mathcal{L}}_{cof}$. Then the binary
relation $\sim^r_\mathcal{L}$ is an equivalence relation on
$\mathcal{M}^{f,\mathcal{L}}_{cof}(X,Y)$.  \ep

Notice that we do not need suppose in the proof of
Proposition~\ref{extright} that the weak factorization system $({\rm
  Cof}\cap \mathcal{W},{\rm Fib})$ is cofibrantly generated. So we do
not need this hypothesis in the proof of
Proposition~\ref{traditionnel1}.

\bpf By Proposition~\ref{extright} applied to the model structure
$({\rm Cof}_\mathcal{L},{\rm
  Fib}_\mathcal{L},\mathcal{W}_\mathcal{L})$, the binary relation
$\sim^r_\mathcal{L}$ coincides with right homotopy for the model
structure $({\rm Cof}_\mathcal{L}, {\rm Fib}_\mathcal{L},
\mathcal{W}_\mathcal{L})$. Since $Y$ is fibrant for the latter model
structure, one deduces that $\sim^r_\mathcal{L}$ is transitive by
\cite{MR99h:55031} Proposition~1.2.5.  \epf

\begin{cor}
  If $\mathcal{L}$ is the class of trivial cofibrations of a left
  Bousfield localization of the model structure of $\mathcal{M}$, then
  the binary relation $\sim^r_\mathcal{L}$ on the set of morphisms
  $\mathcal{M}(X,Y)$ with $X\in\mathcal{M}_{cof}$ and with $Y$ fibrant
  with respect to $\mathcal{L}$ is an equivalence relation.
\end{cor}

\bp \label{quotient2} (dual to \cite{ideeloc} Lemma~3.2) Let $f,g:X
\rightrightarrows Y$ be two morphisms of $\mathcal{M}$. Let
$u:Y\longrightarrow U$ and $v:V\longrightarrow X$ be two other
morphisms of $\mathcal{M}$. If $f\sim^r_{\mathcal{L}} g$, then $u\circ
f \sim^r_{\mathcal{L}} u\circ g$ and $f\circ v\sim^r_{\mathcal{L}}
g\circ v$.  \ep

In other terms, the equivalence relation $\sim_{\mathcal{L}}$ defines
a congruence in the sense of \cite{MR1712872}.

\bpf By considering the opposite of the category $\mathcal{M}$, the
proof is complete using \cite{ideeloc} Lemma~3.2.  \epf

The proof of Proposition~\ref{quotient2} does use the factorization
$(\alpha,\beta)$ and its functoriality. We could avoid using the
functoriality since the morphism $\alpha(\id_Y,\id_Y):Y
\longrightarrow \Path_\mathcal{L}(Y)$ belongs to $\mathcal{L}$ and
since the morphism $\beta(\id_U\p \id_U):\Path_\mathcal{L}(U)
\longrightarrow U \p U$ belongs to $\mathcal{R}$. But anyway, the
proof of Proposition~\ref{quotient2} cannot be adapted to the usual
notion of right homotopy. This is once again a difference between our
notion of right homotopy and the usual one on model category.

Proposition~\ref{quotient2} allows to consider the quotients
$\mathcal{M}/\!\sim_{\mathcal{L}}$ (resp.
$\mathcal{M}_{cof}/\!\sim_{\mathcal{L}}$,
$\mathcal{M}^{f,\mathcal{L}}/\!\sim_{\mathcal{L}}$,
$\mathcal{M}^{f,\mathcal{L}}_{cof}/\!\sim_{\mathcal{L}}$) of the
category $\mathcal{M}$ (resp. $\mathcal{M}_{cof}$,
$\mathcal{M}^{f,\mathcal{L}}$ $\mathcal{M}^{f,\mathcal{L}}_{cof}$) by
the congruence $\sim_{\mathcal{L}}$. By definition, the objects of
$\mathcal{M}/\!\sim_{\mathcal{L}}$ (resp.
$\mathcal{M}_{cof}/\!\sim_{\mathcal{L}}$,
$\mathcal{M}^{f,\mathcal{L}}/\!\sim_{\mathcal{L}}$,
$\mathcal{M}^{f,\mathcal{L}}_{cof}/\!\sim_{\mathcal{L}}$) are the
objects of $\mathcal{M}$ (resp. $\mathcal{M}_{cof}$,
$\mathcal{M}^{f,\mathcal{L}}$, $\mathcal{M}^{f,\mathcal{L}}_{cof}$),
and for any object $X$ and $Y$ of $\mathcal{M}$ (resp.
$\mathcal{M}_{cof}$, $\mathcal{M}^{f,\mathcal{L}}$,
$\mathcal{M}^{f,\mathcal{L}}_{cof}$), one has
$\mathcal{M}/\!\sim_{\mathcal{L}}(X,Y) =
\mathcal{M}(X,Y)/\!\sim_{\mathcal{L}}$ (resp.
$\mathcal{M}_{cof}/\!\sim_{\mathcal{L}}(X,Y) =
\mathcal{M}_{cof}(X,Y)/\!\sim_{\mathcal{L}}$,
$\mathcal{M}^{f,\mathcal{L}}/\!\sim_{\mathcal{L}}(X,Y) =
\mathcal{M}^{f,\mathcal{L}}(X,Y)/\!\sim_{\mathcal{L}}$,
$\mathcal{M}^{f,\mathcal{L}}_{cof}/\!\sim_{\mathcal{L}}(X,Y) =
\mathcal{M}^{f,\mathcal{L}}_{cof}(X,Y)/\!\sim_{\mathcal{L}}$). Let
\beas && [-]_{\mathcal{L}}:
\mathcal{M} \longrightarrow \mathcal{M}/\!\sim_{\mathcal{L}} \\
&& \\
&& [-]_{\mathcal{L}}:
\mathcal{M}_{cof} \longrightarrow \mathcal{M}_{cof}/\!\sim_{\mathcal{L}} \\
&& \\
&&[-]_{\mathcal{L}}:
\mathcal{M}^{f,\mathcal{L}} \longrightarrow \mathcal{M}^{f,\mathcal{L}}/\!\sim_{\mathcal{L}}\\
&& \\
&& [-]_{\mathcal{L}}: \mathcal{M}^{f,\mathcal{L}}_{cof}
\longrightarrow \mathcal{M}^{f,\mathcal{L}}_{cof}/\!\sim_{\mathcal{L}}
\eeas be the canonical functors.

\bp \label{fibrantlike} (dual to \cite{ideeloc} Lemma~3.7) Let $f:X
\longrightarrow Y$ be a morphism of $\mathcal{M}$ belonging to
$\mathcal{L}$. Let us suppose that $X$ is fibrant with respect to
${\mathcal{L}}$. Then there exists $g:Y \longrightarrow X$ such that
$f\circ g \sim_{\mathcal{L}} \id_Y$ and $g \circ f = \id_X$.  \ep

\bpf Let us consider the commutative diagram of $\mathcal{M}$
\[
\xymatrix{
X \ar@{->}[dd]_-{f} \ar@{=}[rr] && X \ar@{->}[dd] \\
&& \\
Y \ar@{->}[rr] \ar@{-->}[rruu]^{g}&& \mathbf{1}.}
\]
Since the left vertical arrow is in $\mathcal{L}$ and since the right
vertical arrow is in $\mathcal{R}$ by hypothesis, there exists a lift
$g :Y \longrightarrow X$. In other terms, $g\circ f=\id_X$. The
diagram of $\mathcal{M}$
\[
\xymatrix{
&&& Y \ar@{->}[d]^-{\alpha(\id_Y,\id_Y)} \\
X\ar@{->}[rrru]^{f} \ar@{->}[rrr]_-{\alpha(\id_Y,\id_Y)\circ f}\ar@{->}[dd]_{f} &&& \Path_{\mathcal{L}}(Y) \ar@{->}[dd]^{\beta(\id_Y,\id_Y)} \\
&&& \\
Y \ar@{-->}[rrruu]^{H}\ar@{->}[rrr]_-{(f\circ g,\id_Y)} &&& Y\p Y}
\]
is commutative since
\[\beta(\id_Y,\id_Y)\circ \alpha(\id_Y,\id_Y)\circ
f=(\id_Y,\id_Y)\circ f = (f, f) = (f\circ g,\id_Y)\circ f.\] Since
$f\in\mathcal{L}$ by hypothesis and since
$\beta(\id_Y,\id_Y)\in\mathcal{R}$, there exists $H:Y \longrightarrow
\Path_{\mathcal{L}}(Y)$ preserving the diagram above commutative.  The
morphism $H$ is by construction a right homotopy from $f\circ g$ to
$\id_Y$ with respect to $\mathcal{L}$.  \epf

\bp \label{petitfibrant} (almost dual to \cite{ideeloc} Theorem~3.9)
One has the isomorphism of categories
$\mathcal{M}^{f,\mathcal{L}}_{cof}/\!\sim_\mathcal{L} \iso
\mathcal{M}^{f,\mathcal{L}}_{cof}[\mathcal{L}^{-1}]$.  In particular,
this means that the category
$\mathcal{M}^{f,\mathcal{L}}_{cof}[\mathcal{L}^{-1}]$ is locally
small.  \ep

The proof of Proposition~\ref{petitfibrant} also shows the isomorphism
of categories
\[\mathcal{M}^{f,\mathcal{L}}/\!\sim_\mathcal{L} \iso
\mathcal{M}^{f,\mathcal{L}}[\mathcal{L}^{-1}].\]

\bpf We know that the pair $(\mathcal{L},\mathcal{R})$ restricts to a
weak factorization system of $\mathcal{M}^{f,\mathcal{L}}_{cof}$. By
considering the opposite category, the proposition is then a
consequence of \cite{ideeloc} Theorem~3.9.  \epf

\bp (Detecting weak equivalences) \label{detection1} A morphism
$f:A\longrightarrow B$ of $\mathcal{M}^{f,\mathcal{L}}_{cof}$ is an
isomorphism of $\mathcal{M}^{f,\mathcal{L}}_{cof}[\mathcal{L}^{-1}]$
if and only if for every object $X$ of
$\mathcal{M}^{f,\mathcal{L}}_{cof}$, the map
\[\mathcal{M}(B,X)/\!\sim_{\mathcal{L}} \longrightarrow \mathcal{M}(A,X)/\!\sim_{\mathcal{L}}\] 
is bijective.  \ep

Note the ``opposite'' characterization
$\mathcal{M}(X,A)/\!\sim_{\mathcal{L}} \longrightarrow
\mathcal{M}(X,B)/\!\sim_{\mathcal{L}}$ also holds. The statement of
the theorem is chosen for having a characterization as close as
possible to the characterization of weak equivalences in a left
Bousfield localization.

\bpf The condition means that the map
\[(\mathcal{M}^{f,\mathcal{L}}_{cof}/\!\sim_{\mathcal{L}})(B,X) 
\longrightarrow (\mathcal{M}^{f,\mathcal{L}}_{cof}/\!\sim_{\mathcal{L}})(A,X)\]
is a bijection. By Yoneda's lemma applied within the locally small
category $\mathcal{M}^{f,\mathcal{L}}_{cof}/\!\sim_{\mathcal{L}}$, the
condition is equivalent to saying that $f:A\longrightarrow B$ is an
isomorphism of
$\mathcal{M}^{f,\mathcal{L}}_{cof}/\!\sim_{\mathcal{L}}$. By
Proposition~\ref{petitfibrant}, the condition is equivalent to saying
that $f:A\longrightarrow B$ is an isomorphism of
$\mathcal{M}^{f,\mathcal{L}}_{cof}[\mathcal{L}^{-1}]$.
\epf

\bd \label{defreplacement} Let $X$ be an object of
$\mathcal{M}_{cof}$. The {\rm fibrant replacement of $X$ with respect
  to $\mathcal{L}$} is the functorial factorization
\[
\xymatrix{
X \ar@{->}[rr]^-{\alpha(f_X)} && R_{\mathcal{L}}(X)
\ar@{->}[rr]^-{\beta(f_X)} && \mathbf{1}}
\] 
of the unique morphism $f_X:X \longrightarrow \mathbf{1}$.  \ed

The mapping $X\mapsto R_{\mathcal{L}}(X)$ is functorial and yields a
functor from $\mathcal{M}_{cof}$ to
$\mathcal{M}^{f,\mathcal{L}}_{cof}$ since the morphism $\alpha(f_X):X
\longrightarrow R_{\mathcal{L}}(X)$ is a cofibration.

\begin{lem} \label{limcompcell} Let $\lambda$ be a limit ordinal. Let
  $X:\lambda \longrightarrow \mathcal{M}$ and $Y:\lambda
  \longrightarrow \mathcal{M}$ be two transfinite sequences. Let $f:X
  \longrightarrow Y$ be a morphism of transfinite sequences such that
  for any $\mu<\lambda$, $f_\mu: X_\mu \longrightarrow Y_\mu$ belongs
  to $\mathcal{L}$. Then $f_\lambda: X_\lambda \longrightarrow
  Y_\lambda$ belongs to $\mathcal{L}$.  Moreover, if for any
  $\mu<\lambda$, $f_\mu: X_\mu \longrightarrow Y_\mu$ belongs to
  $\cell(J\cup K)$, then $f_\lambda: X_\lambda \longrightarrow
  Y_\lambda$ belongs to $\cell(J\cup K)$ as well.
\end{lem}

Lemma~\ref{limcompcell} and Proposition~\ref{inclusionOK} are very
close to \cite{ref_model2} Proposition~12.4.7. The difference is that
we do not suppose here that the underlying model category is cellular.

\bpf Let $T_0=X_\lambda$. Let us consider the unique transfinite
sequence $T:\lambda \longrightarrow \mathcal{M}$ such that one has the
pushout diagram
\[
\xymatrix{
X_\mu \ar@{->}[rr]^{f_\mu}\ar@{->}[dd] && Y_\mu\ar@{->}[dd]\\
&&\\
T_\mu \ar@{->}[rr] && T_{\mu+1}\cocartesien}
\] 
where the left vertical arrow is the composite $X_\mu \longrightarrow
X_\lambda \longrightarrow T_\mu$ for any $\mu<\lambda$. Let $Z$ be an
object of $\mathcal{M}$ and let $\phi:Y_\lambda \longrightarrow Z$ be
a morphism of $\mathcal{M}$. The composite $X_\lambda \longrightarrow
Y_\lambda \longrightarrow Z$ together with the composite $Y_0
\longrightarrow Y_\lambda \longrightarrow Z$ yields with the pushout
diagram above for $\mu=0$ a morphism $T_1 \longrightarrow Z$ since $f$
is a morphism of transfinite sequences.  And by an immediate
transfinite induction, one obtains a morphism $\liminj_\mu T_\mu
\longrightarrow Z$. So one has the isomorphism $\liminj_\mu T_\mu \iso
Y_\lambda$ since the two objects of $\mathcal{M}$ satisfy the same
universal property. Hence the result since the class of morphisms
$\mathcal{L}$ and $\cell(J\cup K)$ are both closed under transfinite
composition.  \epf

\bp \label{inclusionOK} One has $R_\mathcal{L}(\mathcal{L}) \subset
\mathcal{L}$, and even $R_\mathcal{L}(\cell(J\cup K)) \subset
\cell(J\cup K)$.  \ep

\bpf A morphism $f\in\cof(J\cup K)$ is a retract of a morphism $g\in
\cell(J\cup K)$ since $J\cup K$ permits the small object argument. And
the morphism $R_{\cof(J\cup K)}(f)$ is then a retract of the morphism
$R_{\cof(J\cup K)}(g)$. Therefore it suffices to prove that $f\in
\cell(J\cup K)$ implies $R_{\cof(J\cup K)}(f)\in \cell(J\cup K)$. The
functor $R_{\cof(J\cup K)}$ is obtained by a transfinite construction
involving the small object argument. Let $X_0=X$ and $Y_0=Y$ and
$f=f_0$. For any ordinal $\lambda$, let $\overline{Y_\lambda}$ be the
object of $\mathcal{M}$ defined by the following commutative diagram:
\[
\xymatrix{
X_0 \ar@{->}[rr]^{f} \ar@{->}[dd]&& Y_0 \ar@{=}[rr] \ar@{->}[dd] && Y_0 \ar@{->}[dd]\\
&&&&\\
X_\lambda \ar@{->}[rr]&& \overline{Y_\lambda} \ar@{->}[rr]\cocartesien && Y_\lambda}
\]
Let us suppose $f_\lambda:X_\lambda \longrightarrow Y_\lambda$
constructed for some $\lambda\geq 0$ and let us suppose that the
morphism $\overline{Y_\lambda} \longrightarrow Y_\lambda$ is an
element of $\cell(J\cup K)$. The small object argument consists of
considering the sets of commutative squares $\{k \longrightarrow
f_{X_\lambda}, k\in J\cup K\}$ and $\{k \longrightarrow
f_{Y_\lambda},k\in J\cup K\}$ where $f_{X_\lambda}:X_\lambda
\longrightarrow \mathbf{1}$ and $f_{Y_\lambda}:Y_\lambda
\longrightarrow \mathbf{1}$ are the canonical morphisms from
respectively $X_\lambda$ and $Y_\lambda$ to the terminal object of
$\mathcal{M}$. The morphism $f_\lambda$ allows the identification of
$\{k \longrightarrow f_{X_\lambda}, k\in J\cup K\}$ with a subset of
$\{k \longrightarrow f_{Y_\lambda},k\in J\cup K\}$. And the morphism
$f_{\lambda+1}:X_{\lambda+1} \longrightarrow Y_{\lambda+1}$ is
obtained by the diagram (where the notations $\dom(k)$ and $\codom(k)$
mean respectively domain and codomain of $k$):
\[
\xymatrix{
\bigsqcup_{\{k
\longrightarrow f_{X_\lambda}, k\in J\cup K\}} \dom(k) \fr{} \fd{} & X_\lambda \fd{}\fr{}   & Y_\lambda \fd{} \\
\bigsqcup_{\{k\longrightarrow f_{X_\lambda}, k\in J\cup K\}} \codom(k)\fr{}  & X_{\lambda+1} \fr{} \cocartesien & \cocartesien\overline{Y_\lambda} \fd{}\\
& \bigsqcup_{\{k\longrightarrow f_{Y_\lambda}, k\in J\cup K\}\backslash \{k\longrightarrow f_{X_\lambda}, k\in J\cup K\}} \ar@{->}[ru] \dom(k) \fd{} & \cocartesien Y_{\lambda+1} \\
&\bigsqcup_{\{k\longrightarrow f_{Y_\lambda}, k\in J\cup K\}\backslash \{k\longrightarrow f_{X_\lambda}, k\in J\cup K\}}
\ar@{->}[ru] \codom(k) & }
\]
Therefore $f_{\lambda+1}:X_{\lambda+1} \longrightarrow Y_{\lambda+1}$
is an element of $\cell(J\cup K)$. The proof is complete with
Lemma~\ref{limcompcell}.  \epf

Note the same kind of argument as the one of
Proposition~\ref{inclusionOK} leads to the following proposition
(worth being noticed, but useless for the sequel):

\bp One has $\Path_\mathcal{L}(\mathcal{L}) \subset \mathcal{L}$, and
even $\Path_\mathcal{L}(\cell(J\cup K)) \subset \cell(J\cup K)$.  \ep

Proposition~\ref{inclusionOK} will be used in particular in the proof
of Proposition~\ref{prefin1} with the functorial weak factorization
system $({\rm Cof}\cap \mathcal{W},{\rm Fib})$ and in the proof of
Proposition~\ref{pareilfibrant}.

\bp \label{pareilfibrant} The inclusion functor
$\mathcal{M}_{cof}^{f,\mathcal{L}} \subset \mathcal{M}_{cof}$ induces
an equivalence of categories
$\mathcal{M}_{cof}^{f,\mathcal{L}}[\mathcal{L}^{-1}] \simeq
\mathcal{M}_{cof}[\mathcal{L}^{-1}]$. In particular, this implies that
the category $\mathcal{M}_{cof}[\mathcal{L}^{-1}]$ is locally small.
\ep

\bpf Since $R_\mathcal{L}(\mathcal{L})\subset \mathcal{L}$ by
Proposition~\ref{inclusionOK}, there exists a unique functor
$L(R_\mathcal{L})$ making the following diagram commutative:
\[
\xymatrix{
\mathcal{M}_{cof} \ar@{->}[rrr]^-{R_\mathcal{L}} \fd{} &&& \mathcal{M}_{cof}^{f,\mathcal{L}} \fd{} \\
\mathcal{M}_{cof}[\mathcal{L}^{-1}] \ar@{-->}[rrr]^-{L(R_\mathcal{L})} &&& \mathcal{M}_{cof}^{f,\mathcal{L}}[\mathcal{L}^{-1}].
}
\] 
If $i:\mathcal{M}_{cof}^{f,\mathcal{L}} \longrightarrow
\mathcal{M}_{cof}$ is the inclusion functor, then there exists a
unique functor $L(i)$ making the following diagram commutative:
\[
\xymatrix{
\mathcal{M}_{cof}^{f,\mathcal{L}} \ar@{->}[rrr]^-{i} \fd{} &&& \mathcal{M}_{cof} \fd{} \\
\mathcal{M}_{cof}^{f,\mathcal{L}}[\mathcal{L}^{-1}] \ar@{-->}[rrr]^-{L(i)} &&& \mathcal{M}_{cof}[\mathcal{L}^{-1}]. 
}
\] 
There are two natural transformations
$\mu:\id_{\mathcal{M}_{cof}}\Rightarrow i \circ R_{\mathcal{L}}$ and
$\nu:\id_{\mathcal{M}_{cof}^{f,\mathcal{L}}}\Rightarrow
R_{\mathcal{L}} \circ i$ such that for any $X\in \mathcal{M}_{cof}$
and any $Y\in \mathcal{M}_{cof}^{f,\mathcal{L}}$, the morphisms
$\mu(X)$ and $\nu(Y)$ belong to $\mathcal{L}$. So at the level of
localizations, one obtains the isomorphisms of functors
$\id_{\mathcal{M}_{cof}[\mathcal{L}^{-1}]} \iso L(i) \circ
L(R_{\mathcal{L}})$ and
$\id_{\mathcal{M}_{cof}^{f,\mathcal{L}}[\mathcal{L}^{-1}]} \iso
L(R_{\mathcal{L}}) \circ L(i)$. Hence the result.  \epf

\bp \label{incset} Let $(\mathcal{L}',\mathcal{R}')$ be another
cofibrantly generated weak factorization of $\mathcal{M}$ such that
${\rm Cof}\cap \mathcal{W} \subset \mathcal{L}' \subset {\rm Cof}$.
Let us suppose that $\mathcal{L}' \subset\mathcal{L}$. Then the
localization functor $\mathcal{M}_{cof} \longrightarrow
\mathcal{M}_{cof}[\mathcal{L}^{-1}]$ factors uniquely as a composite
\[\mathcal{M}_{cof} \longrightarrow
\mathcal{M}_{cof}[\mathcal{L}'^{-1}] \longrightarrow
\mathcal{M}_{cof}[\mathcal{L}^{-1}].\] 
\ep

\bpf One has $\mathcal{L}' \subset \mathcal{L}$.  \epf

\bp \label{prefin1} The localization functor $L:\mathcal{M}_{cof}
\longrightarrow \mathcal{M}_{cof}[\mathcal{L}^{-1}]$ sends the weak
equivalences of $\mathcal{M}$ between cofibrant objects to
isomorphisms of $\mathcal{M}_{cof}[\mathcal{L}^{-1}]$.  \ep

\bpf The localization functor $L:\mathcal{M}_{cof} \longrightarrow
\mathcal{M}_{cof}[\mathcal{L}^{-1}]$ factors uniquely as a composite
\[\mathcal{M}_{cof} \longrightarrow \mathcal{M}_{cof}[({\rm Cof}\cap
\mathcal{W})^{-1}] \longrightarrow
\mathcal{M}_{cof}[\mathcal{L}^{-1}]\] by Proposition~\ref{incset} and
since ${\rm Cof}\cap \mathcal{W} \subset \mathcal{L}$. By
Proposition~\ref{pareilfibrant} and Proposition~\ref{inclusionOK}
applied to $\mathcal{L} = {\rm Cof}\cap \mathcal{W}$, one has the
equivalence of categories
\[\mathcal{M}_{cof}[({\rm Cof}\cap
\mathcal{W})^{-1}] \simeq \mathcal{M}_{cof}^{f,{\rm Cof}\cap
  \mathcal{W}}[({\rm Cof}\cap \mathcal{W})^{-1}].\] By
Proposition~\ref{petitfibrant} applied to $\mathcal{L}={\rm Cof}\cap
\mathcal{W}$, one has the isomorphism of categories
\[\mathcal{M}_{cof}^{f,{\rm Cof}\cap \mathcal{W}}[({\rm Cof}\cap
\mathcal{W})^{-1}] \iso \mathcal{M}_{cof}^{f,{\rm Cof}\cap
  \mathcal{W}}/\!\sim_{{\rm Cof}\cap \mathcal{W}}.\] Therefore one
obtains the equivalence of categories \[\mathcal{M}_{cof}[({\rm
  Cof}\cap \mathcal{W})^{-1}] \simeq \mathcal{M}_{cof}^{f,{\rm
    Cof}\cap \mathcal{W}}/\!\sim_{{\rm Cof}\cap \mathcal{W}}.\] The
category $\mathcal{M}_{cof}^{f,{\rm Cof}\cap \mathcal{W}}$ is the full
subcategory of cofibrant-fibrant objects of $\mathcal{M}$.  Since
right homotopy with respect to ${\rm Cof}\cap \mathcal{W}$ corresponds
to the usual notion of right homotopy by Proposition~\ref{extright},
one has the equivalence of categories
\[\mathcal{M}_{cof}^{f,{\rm Cof}\cap
\mathcal{W}}/\!\sim_{{\rm Cof}\cap
\mathcal{W}} \simeq \ho(\mathcal{M})\] where 
$\ho(\mathcal{M})=\mathcal{M}[\mathcal{W}^{-1}]$. Hence the result.
\epf

\bp \label{encoredeux} The categories
$\mathcal{M}_{cof}[\mathcal{L}^{-1}]$ and
$\mathcal{M}_{cof}[(\mathcal{W}\cup\mathcal{L})^{-1}]$ are isomorphic.
In particular, this implies that the category
$\mathcal{M}_{cof}[(\mathcal{W}\cup\mathcal{L})^{-1}]$ is locally
small.  \ep

\bpf By Proposition~\ref{prefin1}, there exists a unique functor
\[\mathcal{M}_{cof}[(\mathcal{W}\cup\mathcal{L})^{-1}] \longrightarrow
\mathcal{M}_{cof}[\mathcal{L}^{-1}]\] such that the following diagram
is commutative:
\[
\xymatrix{
\mathcal{M}_{cof} \fd{} \ar@{=}[rr] && \mathcal{M}_{cof} \fd{} \\
\mathcal{M}_{cof}[(\mathcal{W}\cup\mathcal{L})^{-1}] \ar@{-->}[rr] && \mathcal{M}_{cof}[\mathcal{L}^{-1}].
}
\]
Since $\mathcal{L}\subset \mathcal{W}\cup\mathcal{L}$, there exists a
unique functor $\mathcal{M}_{cof}[\mathcal{L}^{-1}] \longrightarrow
\mathcal{M}_{cof}[(\mathcal{W}\cup\mathcal{L})^{-1}]$ such that the
following diagram is commutative:
\[
\xymatrix{
\mathcal{M}_{cof} \fd{} \ar@{=}[rr] && \mathcal{M}_{cof} \fd{} \\
\mathcal{M}_{cof}[\mathcal{L}^{-1}] \ar@{-->}[rr] && \mathcal{M}_{cof}[(\mathcal{W}\cup\mathcal{L})^{-1}].
}
\]
Hence the result.  \epf

In the same way, one can prove the

\bp \label{encoreune} The categories
$\mathcal{M}^{f,\mathcal{L}}_{cof}[\mathcal{L}^{-1}]$ and
$\mathcal{M}^{f,\mathcal{L}}_{cof}[(\mathcal{W}\cup\mathcal{L})^{-1}]$
are isomorphic. In particular, this implies that the category
$\mathcal{M}^{f,\mathcal{L}}_{cof}[(\mathcal{W}\cup\mathcal{L})^{-1}]$
is locally small.  \ep

\begin{nota}
Let 
\[\mathcal{W}_\mathcal{L}=\left\{f:X\longrightarrow Y, \forall Z\in \mathcal{M}_{cof}^{f,\mathcal{L}}, \mathcal{M}(R_\mathcal{L}(Q(Y)),Z)/\!\sim_\mathcal{L} \iso \mathcal{M}(R_\mathcal{L}(Q(X)),Z)/\!\sim_\mathcal{L}\right\}.\] 
\end{nota}

\bp \label{fin1} The inclusion functor $\mathcal{M}_{cof}\subset
\mathcal{M}$ induces an equivalence of categories
$\mathcal{M}_{cof}[(\mathcal{W}\cup \mathcal{L})^{-1}] \simeq
\mathcal{M}[\mathcal{W}_\mathcal{L}^{-1}]$.  In particular, this
implies that the category $\mathcal{M}[\mathcal{W}_\mathcal{L}^{-1}]$
is locally small.  \ep

\bpf Let us consider the composite
\[
\xymatrix{
\mathcal{M} \fr{Q} & \mathcal{M}_{cof} \fr{L} & \mathcal{M}_{cof}[(\mathcal{W}\cup \mathcal{L})^{-1}]}
\] 
where $Q$ is the cofibrant replacement functor of $\mathcal{M}$. By
definition of $\mathcal{W}_\mathcal{L}$ and by
Proposition~\ref{detection1}, by Proposition~\ref{pareilfibrant} and
by Proposition~\ref{encoredeux}, the functor $L\circ Q$ sends the
morphisms of $\mathcal{W}_\mathcal{L}$ to isomorphisms. Thus, there
exists a unique functor $L(Q)$ making the following diagram
commutative:
\[
\xymatrix{
\mathcal{M} \ar@{->}[rrr]^-{Q} \fd{} &&& \mathcal{M}_{cof} \fd{} \\
\mathcal{M}[\mathcal{W}_\mathcal{L}^{-1}] \ar@{-->}[rrr]^-{L(Q)} &&& \mathcal{M}_{cof}[(\mathcal{W}\cup \mathcal{L})^{-1}].}
\] 
Let $i:\mathcal{M}_{cof} \longrightarrow \mathcal{M}$ be the inclusion
functor. Let $f:X\longrightarrow Y\in \mathcal{W}\cup \mathcal{L}$ be
a morphism of $\mathcal{M}_{cof}$. Then $Q(f):Q(X) \longrightarrow
Q(Y)$ is still invertible in $\mathcal{M}_{cof}[(\mathcal{W}\cup
\mathcal{L})^{-1}]$ since the morphism $Q(X)\longrightarrow X$ and
$Q(Y)\longrightarrow Y$ are both weak equivalences of $\mathcal{M}$
between cofibrant objects. So by Proposition~\ref{pareilfibrant}, the
morphism $R_\mathcal{L}(Q(f)):R_\mathcal{L}(Q(X)) \longrightarrow
R_\mathcal{L}(Q(Y))$ is invertible in
$\mathcal{M}^{f,\mathcal{L}}_{cof}[(\mathcal{W} \cup
\mathcal{L})^{-1}] \simeq
\mathcal{M}^{f,\mathcal{L}}_{cof}[\mathcal{L}^{-1}]$. So by
Proposition~\ref{detection1}, one deduces that $f\in
\mathcal{W}_\mathcal{L}$. Thus, there exists a unique functor $L(i)$
making the following diagram commutative:
\[
\xymatrix{
\mathcal{M}_{cof} \ar@{->}[rrr]^-{i} \fd{} &&& \mathcal{M} \fd{} \\
\mathcal{M}_{cof}[(\mathcal{W}\cup \mathcal{L})^{-1}] \ar@{-->}[rrr]^-{L(i)} &&& 
\mathcal{M}[\mathcal{W}_\mathcal{L}^{-1}]. }
\] 
There exist two natural transformations $\mu:Q\circ i \Rightarrow
\id_{\mathcal{M}_{cof}}$ and $\nu:i\circ Q \Rightarrow
\id_{\mathcal{M}}$. If $X$ is an object of $\mathcal{M}_{cof}$, then
$\mu(X)$ is a trivial fibration between cofibrant objects, i.e.
$\mu(X)\in \mathcal{W}$. So one deduces that $L(\mu(X))$ is an
isomorphism of $\mathcal{M}_{cof}[(\mathcal{W}\cup
\mathcal{L})^{-1}]$. Therefore one obtains the isomorphism of functors
$L(Q)\circ L(i) \iso
\id_{\mathcal{M}_{cof}[(\mathcal{W}\cup\mathcal{L})^{-1}]}$. If $Y$ is
an object of ${\mathcal{M}}$, then $\nu(Y):Q(Y)\longrightarrow Y$ is a
trivial fibration of $\mathcal{M}$. Thus, $Q(\nu(Y)):Q(Q(Y))
\longrightarrow Q(Y)$ is a trivial cofibration of $\mathcal{M}$
between cofibrant objects. So $Q(\nu(Y))\in {\rm Cof}\cap\mathcal{W}
\subset \mathcal{L}$. Since $R_\mathcal{L}(\mathcal{L})\subset
\mathcal{L}$, one deduces that $R_\mathcal{L}(Q(\nu(Y)))$ is an
isomorphism of $\mathcal{M}^{f,\mathcal{L}}_{cof}[\mathcal{L}^{-1}]$.
Again by Proposition~\ref{detection1}, one deduces that $\nu(Y)\in
\mathcal{W}_\mathcal{L}$ and one obtains the isomorphism of functors
$L(i)\circ L(Q) \iso \id_{\mathcal{M}[\mathcal{W}_\mathcal{L}^{-1}]}$.
The proof is complete.  \epf

\bth (Whitehead's theorem for the localization of a model category
with respect to a weak factorization system)
\label{fin2}
The inclusion functors $\mathcal{M}_{cof}^{f,\mathcal{L}} \subset
\mathcal{M}_{cof} \subset \mathcal{M}$ induce the equivalences of
categories
\[ 
\mathcal{M}_{cof}^{f,\mathcal{L}}/\!\sim_\mathcal{L} \simeq 
\mathcal{M}_{cof}[(\mathcal{W}\cup \mathcal{L})^{-1}] \simeq 
\mathcal{M}[\mathcal{W}_\mathcal{L}^{-1}]. \] 
The functor $\mathcal{M}[\mathcal{W}_\mathcal{L}^{-1}]
\longrightarrow\mathcal{M}_{cof}^{f,\mathcal{L}}/\!\sim_\mathcal{L}$ is given by the cofibrant-fibrant
w.r.t. $\mathcal{L}$ functor $R_\mathcal{L} \circ Q : \mathcal{M}
\longrightarrow \mathcal{M}_{cof}^{f,\mathcal{L}}$. Moreover, the
localization functor $\mathcal{M} \longrightarrow
\mathcal{M}[\mathcal{W}_\mathcal{L}^{-1}]$ factors uniquely as a
composite
\[\mathcal{M} \longrightarrow \mathcal{M}[\mathcal{W}^{-1}] \longrightarrow \mathcal{M}[\mathcal{W}_\mathcal{L}^{-1}].\]
\eth

The equivalence of categories $\mathcal{M}_{cof}[(\mathcal{W}\cup
\mathcal{L})^{-1}] \simeq \mathcal{M}[\mathcal{W}_\mathcal{L}^{-1}]$
shows that up to weak equivalence and up to the 2-out-of-3 axiom, a
morphism of $\mathcal{W}_\mathcal{L}$ is a morphism of
$\mathcal{W}\cup \mathcal{L}$ between cofibrant objects of
$\mathcal{M}$. This means that the class of morphisms
$\mathcal{W}_\mathcal{L}$ is not too big.

\bpf The equivalence of categories
$\mathcal{M}[\mathcal{W}_\mathcal{L}^{-1}] \simeq
\mathcal{M}_{cof}[(\mathcal{W}\cup \mathcal{L})^{-1}]$ is given by
Proposition~\ref{fin1}. By Proposition~\ref{fin1}, the functor
$\mathcal{M}[\mathcal{W}_\mathcal{L}^{-1}] \longrightarrow
\mathcal{M}_{cof}[(\mathcal{W}\cup \mathcal{L})^{-1}]$ is induced by
the cofibrant replacement functor $Q:\mathcal{M} \longrightarrow
\mathcal{M}_{cof}$. Therefore every morphism of $\mathcal{W}$ is in
$\mathcal{W}_\mathcal{L}$. At last, one has
\begin{align*}
& \mathcal{M}_{cof}[(\mathcal{W}\cup\mathcal{L})^{-1}] & \\
& \iso \mathcal{M}_{cof}[\mathcal{L}^{-1}] & \hbox{ by Proposition~\ref{encoredeux}}\\
& \simeq \mathcal{M}^{f,\mathcal{L}}_{cof}[\mathcal{L}^{-1}] & \hbox{ by Proposition~\ref{pareilfibrant}}\\
& \simeq  \mathcal{M}_{cof}^{f,\mathcal{L}}/\!\sim_\mathcal{L} & \hbox{ by Proposition~\ref{petitfibrant}.}
\end{align*}
\epf

Theorem~\ref{fin2} says that the category
$\mathcal{M}[\mathcal{W}_\mathcal{L}^{-1}]$ inverts all weak
equivalences of $\mathcal{M}$ and all morphisms of $\mathcal{L}$ with
cofibrant domains. We do not know if all morphisms of $\mathcal{L}$
(and not only the ones with cofibrant domain) are inverted in the
category $\mathcal{M}[\mathcal{W}_\mathcal{L}^{-1}]$. But there is a
kind of reciprocal statement:

\bp Let us suppose that $\mathcal{M}$ is left proper. Let $f$ be a
cofibration of $\mathcal{M}$ such that $Q(f)\in \mathcal{L}$. Then
$f\in \mathcal{L}$.  \ep

\bpf Let $p\in\mathcal{R}$. Since ${\rm Cof}\cap\mathcal{W}\subset
\mathcal{L}$, the morphism $p$ is a fibration of $\mathcal{M}$. By
hypothesis, $p$ satisfies the RLP with respect to $Q(f)$. Since $f$ is
a cofibration and by \cite{ref_model2} Proposition~13.2.1, one deduces
that $p$ satisfies the RLP with respect to $f$. So $f\in\mathcal{L}$.
\epf

Before treating the case of T-homotopy equivalences in
Section~\ref{hc}, let us give some examples of the situation explored
in this section.

\subsubsection*{Example~1}

Let $\mathcal{M}$ be a cofibrantly generated model category with set
of generating cofibrations $I$ and with set of generating trivial
cofibrations $J$ and with class of weak equivalences $\mathcal{W}$.
Let $(\mathcal{L},\mathcal{R})=(\cof(J),\inj(J))$. The main theorem
gives the equivalences of categories
\[
\mathcal{M}_{cof}^{fib}/\!\sim_{\cof(J)}
\simeq \mathcal{M}_{cof}[\cof(J)^{-1}]
\simeq 
\mathcal{M}[\mathcal{W}_{\cof(J)}^{-1}] 
\] 
where $\mathcal{M}_{cof}^{fib}$ is the full subcategory of
cofibrant-fibrant objects. One can directly check that
$\mathcal{W}_{\cof(J)}^{-1}=\mathcal{W}$. This is not surprising since
the category $\mathcal{M}_{cof}^{fib}/\!\sim_{\cof(J)}$ is the
category of cofibrant-fibrant objects of $\mathcal{M}$ up to homotopy.

\subsubsection*{Example~2}

Let $\mathcal{M}$ be a cofibrantly generated model category with set
of generating cofibrations $I$ and with set of generating trivial
cofibrations $J$ and with class of weak equivalences $\mathcal{W}$.
Then one can consider the model structure $(\All,\All,\Iso)$ where all
morphisms are a cofibration and a fibration and where the weak
equivalences are the isomorphisms. Indeed, one has $(\Iso,\All) =
(\cof(\id_\varnothing),\inj(\id_\varnothing))$. The main theorem
applied with the latter model structure and with the weak
factorization system $(\mathcal{L},\mathcal{R})=(\cof(J),\inj(J))$
gives the equivalences of categories
\[
\mathcal{M}^{fib}/\!\sim_{\cof(J)}
\simeq \mathcal{M}[\cof(J)^{-1}]\simeq 
\mathcal{M}[\Iso_{\cof(J)}^{-1}] 
\] 
where $\mathcal{M}^{fib}$ is the full subcategory of fibrant objects,
where $\sim_{\cof(J)}$ is a congruence on the morphisms of the full
subcategory of fibrant objects of $\mathcal{M}$. The functor
$R_{\cof(J)}$ is a fibrant replacement functor of $\mathcal{M}$ and
the functor $Q$ is a cofibrant replacement functor of the model
structure $(\All,\All,\Iso)$. That is one can suppose that
$Q=\id_\mathcal{M}$.

\subsubsection*{Example~3}

Let $\mathcal{M}$ be a combinatorial model category (in the sense of
Jeff Smith), that is a cofibrantly generated model category such that
the underlying category is locally presentable \cite{MR95j:18001}. Let
$J$ be the set of generating trivial cofibrations. Then for any set
$K$ of morphisms of $\mathcal{M}$, the pair $(\cof(J\cup K),\inj(J\cup
K))$ is a cofibrantly generated weak factorization system by
\cite{MR1780498} Proposition~1.3. Then the main theorem of this
section applies. One obtains the equivalences of categories
\[\mathcal{M}_{cof}^{f,{\cof(J\cup K)}}/\!\sim_{\cof(J\cup K)} 
\simeq \mathcal{M}_{cof}[(\mathcal{W}\cup {\cof(J\cup K)})^{-1}]
\simeq \mathcal{M}[\mathcal{W}_{\cof(J\cup K)}^{-1}] .\] Assume
Vop\v{e}nka's principle (\cite{MR95j:18001} chapter~6). If
$\mathcal{M}$ is left proper, then the left Bousfield localization
$L_K \mathcal{M}$ of the model category $\mathcal{M}$ with respect to
the set of morphisms $K$ exists by a theorem of Jeff Smith proved in
\cite{MR1780498} Theorem~1.7 and in \cite{rotho} Theorem~2.2. The
category $\mathcal{M}[\mathcal{W}_{\cof(J\cup K)}^{-1}]$ is not
necessarily equivalent to the Quillen homotopy category
$\ho(L_K\mathcal{M})$ of this Bousfield localization. But all
morphisms inverted by $\mathcal{M}[\mathcal{W}_{\cof(J\cup K)}^{-1}]$
are inverted by the Bousfield localization. In other terms, the
functor $\mathcal{M} \longrightarrow \ho(L_K\mathcal{M})$ factors
uniquely as a composite $\mathcal{M} \longrightarrow
\mathcal{M}[\mathcal{W}_{\cof(J\cup K)}^{-1}] \longrightarrow
\ho(L_K\mathcal{M})$.

\subsubsection*{Example~4}

Let $\mathcal{M}$ be a model category with model structure denoted by
$({\rm Cof},{\rm Fib},\mathcal{W})$ for respectively the class of
cofibrations, of fibrations and of weak equivalences such that the
weak factorization system $({\rm Cof}\cap \mathcal{W},{\rm Fib})$ is
cofibrantly generated : the set of generating trivial cofibrations is
denoted by $J$. Let $(\mathcal{L},\mathcal{R})$ be a cofibrantly
generated weak factorization system such that ${\rm Cof}\cap
\mathcal{W}\subset \mathcal{L} \subset {\rm Cof}$. Let us suppose that
the left Bousfield localization $L_\mathcal{L}\mathcal{M}$ of
$\mathcal{M}$ with respect to $\mathcal{L}$ exists and let us suppose
that $\mathcal{L}$ is the class of trivial cofibrations of this
Bousfield localization. One obtains the equivalences of categories
\[\mathcal{M}_{cof}^{f,\mathcal{L}}/\!\sim_\mathcal{L} \simeq 
\mathcal{M}_{cof}[(\mathcal{W}\cup \mathcal{L})^{-1}]\simeq 
\mathcal{M}[\mathcal{W}_\mathcal{L}^{-1}].\] 
The category $\mathcal{M}_{cof}^{f,\mathcal{L}}$ is the full
subcategory of $\mathcal{M}$ containing the cofibrant-fibrant object
of $L_\mathcal{L}\mathcal{M}$. The congruence $\sim_\mathcal{L}$ is
the usual notion of homotopy in $\mathcal{M}$ (\cite{ref_model2}
Proposition~3.5.3). Then the category 
$\mathcal{M}_{cof}^{f,\mathcal{L}}/\!\sim_\mathcal{L}$ is equivalent
to the Quillen homotopy category $\ho(L_\mathcal{L}\mathcal{M})$.

\section{Application : homotopy continuous flow and Whitehead's theorem}
\label{hc}

The category $\top$ of \textit{compactly generated topological spaces}
(i.e. of weak Hausdorff $k$-spaces) is complete, cocomplete and
cartesian closed (more details for this kind of topological spaces in
\cite{MR90k:54001,MR2000h:55002}, the appendix of \cite{Ref_wH} and
also the preliminaries of \cite{model3}). For the sequel, all
topological spaces will be supposed to be compactly generated. A
\textit{compact space} is always Hausdorff. The category $\top$ is
equipped with the unique model structure having the \textit{weak
  homotopy equivalences} as weak equivalences and having the
\textit{Serre fibrations}~\footnote{that is a continuous map having
  the RLP with respect to the inclusion $\mathbf{D}^n\p 0\subset
  \mathbf{D}^n\p [0,1]$ for all $n\geq 0$ where $\mathbf{D}^n$ is the
  $n$-dimensional disk.} as fibrations \cite{MR99h:55031}.

As already described in the introduction, the time flow of a higher
dimensional automaton is encoded in an object called a \textit{flow}.
The category $\dtop$ is equipped with the unique model structure such
that \cite{model3}:
\begin{itemize}
\item The weak equivalences are the \textit{weak S-homotopy equivalences}, 
i.e. the morphisms of flows $f:X\longrightarrow Y$ such that
$f^0:X^0\longrightarrow Y^0$ is a bijection and such that $\P f:\P
X\longrightarrow \P Y$ is a weak homotopy equivalence. 
\item The fibrations are the morphisms of flows $f:X\longrightarrow Y$
  such that $\P f:\P X \longrightarrow \P Y$ is a Serre fibration.
\end{itemize}
This model structure is cofibrantly generated. The set of generating
cofibrations is the set $I^{gl}_+ = I^{gl}\cup \{R,C\}$ with
\[I^{gl}=\{\glob(\mathbf{S}^{n-1})\subset \glob(\mathbf{D}^{n}), n\geq
0\}\] where $\mathbf{D}^{n}$ is the $n$-dimensional disk, where
$\mathbf{S}^{n-1}$ is the $(n-1)$-dimensional sphere, where $R$ and
$C$ are the set maps $R:\{0,1\} \longrightarrow \{0\}$ and
$C:\varnothing \longrightarrow \{0\}$ and where for any topological
space $Z$, the flow $\glob(Z)$ is the flow defined by
$\glob(Z)^0=\{\widehat{0},\widehat{1}\}$, $\P \glob(Z)=Z$,
$s=\widehat{0}$ and $t=\widehat{1}$, and a trivial composition law.
The set of generating trivial cofibrations is
\[J^{gl}=\{\glob(\mathbf{D}^{n}\p\{0\})\subset
\glob(\mathbf{D}^{n}\p [0,1]), n\geq 0\}.\]

The weak S-homotopy model structure of $\dtop$ has some similarity
with the model structure on the category of small simplicial
categories (with identities !) constructed in \cite{math.AT/0406507}.
The weak equivalences (resp. the fibrations) of the latter look like
the weak equivalences (resp. the fibrations) of the model structure of
$\dtop$ with an additional condition. The weak S-homotopy model
structure of $\dtop$ has also some similarity with the model structure
on the category of small simplicial categories (with identities again
!) on a \textit{fixed} set of objets $O$ constructed in
\cite{MR579087}. For the latter, the set maps $R:\{0,1\}
\longrightarrow \{0\}$ and $C:\varnothing \longrightarrow \{0\}$ are
not used since the set of objects is fixed.

\bd 
A flow $X$ is {\rm loopless} if for any $\alpha\in X^0$, the space
$\P_{\alpha,\alpha}X$ is empty.
\ed

Recall that a flow is a small category without identities morphisms
enriched over a category of topological spaces.  So the preceding
definition is meaningful.

A poset $(P,\leq)$ can be identified with a loopless flow having $P$
as set of states and such that there exists a non-constant execution
path from $x$ to $y$ if and only if $x < y$. The corresponding flow is
still denoted by $P$. This defines a functor from the full subcategory
of posets whose morphisms are the strictly increasing maps to the full
subcategory of loopless flows. The category of finite bounded posets
is essentially small. Let us choose a small subcategory of
representatives.

\bd \cite{1eme} \label{genT} Let $\mathcal{T}$ be the set of
cofibrations $Q(f):Q(P_1) \longrightarrow Q(P_2)$ such that
$f:P_1\longrightarrow P_2$ is a morphism of posets satisfying the
following conditions:
\begin{enumerate}
\item The posets $P_1$ and $P_2$ are finite and bounded. 
\item The morphism of posets $f:P_1 \longrightarrow P_2$ is
  one-to-one; in particular, if $x$ and $y$ are two elements of $P_1$
  with $x < y$, then $f(x) < f(y)$.
\item One has $f(\min P_1) = \min P_2 = \widehat{0}$ and $f(\max P_1)
  = \max P_2 = \widehat{1}$.
\item The posets $P_1$ and $P_2$ are objects of the chosen small
  subcategory of representatives of the category of finite bounded
  posets.
\end{enumerate}
The set $\T$ is called the set of {\rm generating T-homotopy
  equivalences}.  \ed

The set $\T$ is introduced in \cite{1eme} for modelling T-homotopy as
a refinement of observation. By now, this is the best known definition
of T-homotopy.

\bd A flow $X$ is {\rm homotopy continuous} if the unique morphism of
flows $f_X:X \longrightarrow \mathbf{1}$ belongs to
$\inj(J^{gl}\cup\T)$.  \ed

Notice that $\inj(J^{gl}\cup\T)=\inj(\T)$ because all flows are
fibrant for the weak S-homotopy model structure of $\dtop$.

Let $X$ be a homotopy continuous flow. Then, for instance, consider
the unique morphism $Q(f):Q(\{\widehat{0}<\widehat{1}\})
\longrightarrow Q(\{\widehat{0}<A<\widehat{1})\}$ of $\mathcal{T}$.
For any commutative square
\[
\xymatrix{
Q(\{\widehat{0}<\widehat{1}\}) \ar@{->}[rr]^{\phi}\ar@{->}[dd] && X
 \\ &&\\ Q(\{\widehat{0}<A<\widehat{1}\})
\ar@{-->}[rruu]^{k} &&}
\]
there exists $k:Q(\{\widehat{0}<A<\widehat{1}\})\longrightarrow X$
making the triangle commutative. Therefore, the existence of $k$
ensures that any directed segment of $X$ can always be divided up to
S-homotopy. Roughly speaking, a flow $X$ is homotopy continuous if it
is indefinitely divisible up to S-homotopy.

\bp The pair $(\cof(J^{gl}\cup \T),\inj(J^{gl}\cup \T))$ is a
cofibrantly generated weak factorization system. Moreover, it
satisfies the conditions of Section~1, that is:
\[\cof(J^{gl}) \subset \cof(J^{gl}\cup \T) \subset \cof(I^{gl}_+).\] 
\ep

\bpf For every $g\in J^{gl}\cup \T$, the continuous map $\P g$ is a
closed inclusion of topological spaces. Therefore by \cite{model3}
Proposition~11.5 and by \cite{MR99h:55031} Theorem~2.1.14, the small
object argument applies.  \epf

\begin{nota} 
\beas
&& \mathcal{S}_\T=\mathcal{S}_{\cof(J^{gl}\cup \T)}, \\
&& R_\T=R_{\cof(J^{gl}\cup \T)}, \\
&& \sim_\T=\sim_{\cof(J^{gl}\cup
\T)}, \\
&& \dtop_{cof}^{f,\T}= \dtop_{cof}^{f,\cof(J^{gl}\cup \T)}.
\eeas
\end{nota}

We can now apply Theorem~\ref{fin2} to obtain the theorem:

\begin{thm} \label{fin3} The inclusion functors
  $\dtop_{cof}^{f,\T}\subset \dtop_{cof} \subset \dtop$ induce the
  equivalences of categories
\[\dtop_{cof}^{f,\T}/\!\sim_\T \simeq \dtop_{cof}[(\mathcal{S}\cup \cof(J^{gl}\cup T))^{-1}] \simeq \dtop[\mathcal{S}_\T^{-1}].\]
\end{thm}

It remains to check the invariance of the underlying homotopy type and
of the branching and merging homology theories:

\begin{thm}\label{fin4}
  A morphism of $\mathcal{S}_\T$ preserves the underlying homotopy
  type and the branching and merging homology theories.
\end{thm}

\bpf It has been already noticed above that up to weak S-homotopy and
up to the 2-out-of-3 axiom, a morphism of $\mathcal{S}_\T$ is a
morphism of $\mathcal{S}\cup \cof(J^{gl}\cup \T)$ between cofibrant
flows. Formally, let $f:A \longrightarrow B$ be an element of
$\mathcal{S}_\T$. Then consider the commutative diagram:
\[
\xymatrix{
Q(A) \ar@{->}[dd]_{Q(f)} \ar@{->}[rr]^{\simeq} && A \ar@{->}[dd]_{f}\\
&& \\
Q(B) \ar@{->}[rr] ^{\simeq} && B}
\]
The morphism $Q(f)$ is an isomorphism of $\dtop_{cof}[(\mathcal{S}\cup
\cof(J^{gl}\cup T))^{-1}]$. Therefore it preserves the underlying
homotopy type and the branching and merging homology theories because
any morphism of $\mathcal{S}$ preserves these invariants by
\cite{model2} Proposition~VII.2.5 and by \cite{exbranch} Corollary~6.5
and Corollary~A.11, and because any morphism of $\cof(J^{gl}\cup \T)$
preserves these invariants by \cite{1eme} Theorem~5.2. The morphisms
$Q(A)\longrightarrow A$ and $Q(B)\longrightarrow B$ are weak
S-homotopy equivalences. So both preserve the underlying homotopy type
\cite{model2} Proposition~VII.2.5 and the branching and merging
homology theories \cite{exbranch} Corollary~6.5 and Corollary~A.11.
Hence the result.  \epf

\end{document}